\documentclass[12pt, reqno]{amsart}
 \usepackage{amsmath, amsthm, amscd, amsfonts, amssymb, graphicx, color, float}
\usepackage[bookmarksnumbered, colorlinks, plainpages]{hyperref}
\input{mathrsfs.sty}
\hypersetup{colorlinks=true,linkcolor=red, anchorcolor=green, citecolor=cyan, urlcolor=red, filecolor=magenta, pdftoolbar=true}
\usepackage{soul,cancel}
\usepackage{tikzsymbols}
\usepackage{textcomp}

\theoremstyle{definition}

\theoremstyle{remark}

\numberwithin{equation}{section}

\newcommand{\K}{\mathbb{K}}
\newcommand{\R}{\mathbb{R}}
\newcommand{\C}{\mathbb{C}}
\newcommand{\N}{\mathbb{N}}
\newcommand{\al}{\alpha}

\begin{document}

\title[Dichotomy between finite and infinite dimensional spaces]{Dichotomy between operators acting on finite and infinite dimensional Hilbert spaces}


\author[Bernal]{L. Bernal-Gonz\'alez}
\address[L. Bernal-Gonz\'alez]{\mbox{}\newline \indent Departamento de An\'{a}lisis Matem\'{a}tico, \newline \indent Facultad de Matem\'{a}ticas, \newline \indent
	Instituto de Matem\'aticas Antonio de Castro Brzezicki, \newline \indent
	Universidad de Sevilla, \newline \indent
	Avenida Reina Mercedes, Sevilla, 41080  (Spain).}
\email{lbernal@us.es}


\author[Moslehian]{M.S.~Moslehian}
\address[M.S.~Moslehian]{\mbox{}\newline \indent Department of Pure Mathematics, \newline \indent Center of Excellence in
	Analysis on Algebraic Structures (CEAAS), \newline \indent Ferdowsi University of
	Mashhad, \newline \indent P.O. Box 1159, Mashhad 91775, Iran}
\email{moslehian@um.ac.ir; moslehian@yahoo.com}


\author[Seoane]{J.B.~Seoane-Sep\'{u}lveda}
\address[J.B.~Seoane-Sep\'{u}lveda]{\mbox{}\newline \indent Instituto de Matem\'atica Interdisciplinar (IMI)\newline \indent Departamento de An\'{a}lisis Matem\'{a}tico y Matem\'atica Aplicada \newline \indent
	Facultad de Ciencias Matem\'{a}ticas\newline \indent
	Plaza de Ciencias 3 \newline \indent
	Universidad Complutense de Madrid \newline \indent
	Madrid, 28040, Spain.}
\email{jseoane@mat.ucm.es}




\thanks{The first author was supported by the Plan Andaluz de Investigaci%
\'{o}n de la Junta de Andaluc\'{\i}a FQM-127 Grant P20-00637 and by MICINN
Grant PGC2018-098474-B-C21. The third author was supported by Grant
PGC2018-097286-B-I00.}
\maketitle

\begin{abstract}
In this expository article, we give several examples showing how drastically different can be the behavior of operators acting on finite versus infinite dimensional Hilbert spaces. This essay is written as in such a friendly-reader to show that the situation in the infinite dimensional setting is trickier than the finite one.
\end{abstract}

\maketitle
The notion of Hilbert space is a generalization of that of the Euclidean space $\mathbb{R}^2$, that is a vector space equipped with a scalar product $\langle\cdot,\cdot\rangle$. The norm is defined by
$$\|x\|=\langle x,x\rangle^{1/2}.$$
This notion was first introduced by David Hilbert in the setting of integral equations and named by others after him.  A Hilbert space is a Banach space  (i.e., a vector space equipped with a complete norm $\|\cdot\|$), whose norm satisfies the parallelogram law  $\|x+y\|^2+\|x-y\|^2=2\|x\|^2+2\|y\|^2$. The space $C[0,1]$ of continuous linear functions on the interval $[0,1]$ endowed with the sup-norm $\|f\|=\sup\{|f(t)|: t\in[0,1]\}$  is a Banach space whose norm cannot be deduced from an inner product space since it does not satisfy the parallelogram law for $f(t)=1$ and $g(t)=t$.

Finite dimensional Hilbert spaces are isomorphic to the space $\mathbb{C}^n$ endowed with the inner product $\langle(x_i), (y_i)\rangle=\sum_{i=1}^nx_i\overline{y_i}$.

The infinite dimensional analogue of $\mathbb{C}^n$ is the (separable) Hilbert space $\ell_2=\ell_2(\mathbb{N})$ of all complex sequences $(x_n)$ satisfying $\sum_{n=1}^\infty|x_n|^2<\infty$ under pointwise operations $\alpha (x_n)+(y_n)=(\alpha x_n+y_n)$ and the inner product $$\langle(x_n), (y_n)\rangle=\sum_{n=1}^\infty x_n\overline{y_n}.$$
The standard orthonormal basis $\{e_j: j=1, 2, \ldots\}$ of $\ell_2$ is the direct analogue of the one of $\mathbb{C}^n$. Similarly, one can impose a Hilbert space structure on the linear space $\ell_2(\mathbb{Z})$ consisting of all two-sided sequences of the form $(\ldots, x_{-2},x_{-1},x_0,x_1, x_2,\ldots)$ such that $\sum_{n=-\infty}^\infty|x_n|^2<\infty$.

Let $(\mathscr{H}, \langle\cdot,\cdot\rangle)$ be a Hilbert space.  By $\mathbb{B}(\mathscr{H})$ we denote the algebra of all continuous linear operators on $\mathscr{H}$ equipped with the pointwise-defined operations of addition and multiplication by scalars, while the multiplication is defined as the composition of operators. A linear operator $A: \mathscr{H}\to \mathscr{H}$ is called bounded if $\|Ax\|\leq M\|x\|$ for some $M\geq 0$ and all $x\in\mathscr{H}$; if this is the case, $\|A\|:=\sup\{\|Ax\|: \|x\|=1\}<\infty$ is called the operator norm. The continuity of a linear operator is equivalent to its boundedness in virtue of $\|Ax-Ay\|\leq \|A\|\,\|x-y\|$.
For every operator $A\in \mathbb{B}(\mathscr{H})$, there exists a unique operator $A^*\in \mathbb{B}(\mathscr{H})$, called the adjoint operator of $A$, such that $\langle Ax, y\rangle=\langle x, A^*y\rangle$ for all $ x,y \in{\mathscr H}$. Throughout the paper, a capital letter means a continuous linear operator in ${\mathbb B}({\mathscr H})$, in particular, $I$ denotes the identity operator. When a capital letter denotes a matrix, we explicitly state it.

The set $\sigma(A):=\{\lambda\in\mathbb{C}: A-\lambda I {\rm ~is~not~invertible~in~} \mathbb{B}(\mathscr{H})\}$ is called the spectrum of $A$, which is both nonempty and compact. The numerical range of $A$ is defined and denoted by $W(A)=\{\langle Ax,x\rangle:\|x\|=1\}$.

An operator $A$ is called normal if $A^*A=AA^*$. It is self-adjoint (or hermitian) if $A^*=A$, or equivalently $W(A)\subseteq \mathbb{R}$. It is said to be positive (positive semidefinite) if $W(A)\subseteq [0,\infty)$; the set of all positive semi-definite operators is denoted by $\mathbb{B}(\mathscr{H})_+$. An operator $A$ is idempotent if $A^2=A$. An orthogonal projection is a self-adjoint idempotent.\\
The L\"owner order on the set $\mathbb{B}(\mathscr{H})_h$ of self-adjoint operators is defined by  $A\leq B \iff B-A\in \mathbb{B}(\mathscr{H})_+$.


There are many assertions in (finite dimensional) linear algebra that do not hold in an infinite dimensional Hilbert space; even less is true for general Banach spaces than Hilbert spaces. \\

First of all, let us explain that by the dimension of a linear space (in the algebraic sense) we understand the cardinality of any of its linear (or Hamel) bases, i.e., maximal linearly independent sets. In Banach spaces of finite dimension, such as $\mathbb{C}^n$, the closed unit ball is compact, all subspaces are (topologically) closed, and all norms on the space are equivalent. None of these statements is true anymore in Banach spaces of infinite dimension, such as $\ell_2$. \\

If $\mathscr{H}$ is a vector space endowed with a scalar product and is finite dimensional, then $\mathscr{H}$ au\-to\-ma\-ti\-cal\-ly becomes a Hilbert space.
But if $\mathscr{H}$ is infinite-dimensional, this is not always true. In fact, if ${\rm dim} (\mathscr{H}) = \aleph_0$ (the cardinality of $\N$), then the
Baire category theorem prevents $\mathscr{H}$ to be Hilbert. \\

A Hilbert space $\mathscr{H}$, besides Hamel bases, also possesses the so-called Hilbert (or orthonormal) bases, that is, maximal families of orthogonal norm-one vectors (two vectors are orthogonal if their inner product is zero). If the dimension of $\mathscr{H}$ is finite, then the Gram--Schmidt process allows us to produce a Hilbert basis from a linear basis, and the cardinalities of these bases are the same. If the dimension of $\mathscr{H}$ is infinite, then the cardinality of a Hilbert basis for $\mathscr{H}$ is strictly smaller than the cardinality of a linear basis for $\mathscr{H}$; see \cite{HAL}.

If $A$ is a continuous linear operator on $\ell^2$ or $\mathbb{C}^n$, then it admits a matrix representation, i.e. an infinite (resp., finite) matrix whose $(i,j)$-entry is $\langle Ae_j,e_i\rangle$ for all pairs $i, j$, and the action of $A$ is described by the usual matrix product (evidently, a change of orthonormal basis results in a different matrix representation, and each can be endowed with some norm; see \cite{FRR} for a study of variation of matrix norms as the basis varies).  The converse is true for $\mathbb{M}_n$ in the sense that an arbitrary matrix $A\in \mathbb{M}_n$ corresponds to the linear mapping on $\mathbb{C}^n$ defined by $[z_1, \ldots, z_n]^t\mapsto A[z_1, \ldots, z_n]^t$ via a matrix product. A similar assertion is not valid for for infinite matrices: not any matrix corresponds to a continuous linear operator. In principle, all information about an operator acting on a finite dimensional Hilbert space can be systematically obtained from its matrix representation; the latter in the infinite dimensional case is not useful. \\
As Halmos indicated \cite[Chapter 5]{HAL}, if $\sum_i\sum_j|\lambda_{ij}|^2<\infty$, then there is an operator (matrix, resp.) $A \in\mathbb{B}(\ell^2)$ such that $\lambda_{ij}=\langle Ae_j,e_i\rangle$. Of course, this condition is not necessary. For example, it is not satisfied even in the simplest case of the identity operator.

Thus, we can identify $\mathbb{B}(\mathbb{C}^n)$ with the space $\mathbb{M}_n$ of all $n\times n$ complex matrices in the canonical way.  In this case, if $A= [a_{ij}] \in \mathbb{M}_n$, then $A^* = [\overline{a_{ji}}]$. In addition, $\sigma(A)$ is exactly the set of eigenvalues of $A$, since $A$ is invertible if and only if it is one-to-one.

Now we present several examples to demonstrate some differences between the properties of operators on finite dimensional Hilbert spaces and those on infinite dimensional ones. It is worthy to say that there are several tricks with matrices, in particular $2\times 2$ ones, which help researchers to establish results concerning operators that could not be treated easily; see e.g. \cite{BHA, MOS2}.

\begin{itemize}
\item A linear operator $A\in \mathbb{B}(\mathbb{C}^n)$ is injective (one-to-one) if and only if it is surjective. This is not the case for linear operators on infinite dimensional Hilbert spaces. For example, the right (unilateral) shift operator $A: \ell_2 \to \ell_2$ defined by \[A(x_1, x_2, \ldots)=(0, x_1, x_2, \ldots)\] is injective but not surjective. In addition, \[A^*(x_1, x_2, \ldots)=(x_2, x_3, \ldots),\] which is called left (backward) shift operator, is surjective but not injective.

From another point of view, we can describe the situation above by stating that a matrix $A$ is an isometry (i.e. $\|Ax||=\|x\|$ for all $x$) if and only if it is unitary.  In the framework of infinite dimensional Hilbert spaces this is not valid, since the right shift operator $A$ is an isometry ($A^*A=I$) but not unitary ($AA^*\neq I$); see also \cite{GO}.

Still, there is another direction to look at this from: We observe that the right shift operator $A$ has a left inverse but not a right inverse whilst a square matrix having a left inverse will automatically have a right inverse.\\

\item Every matrix has an eigenvalue while the right shift operator $A$ has no eigenvalues since $Ax=\lambda x$ implies that $x=0$. This shows that the spectrum of an operator may have no eigenvalue but still is nonempty. It is worthy to mention that the lack of eigenvalues for normal operators is replaced by the spectral theorem.\\

\item By the rank-nullity theorem, $\dim\ker(A)=\dim\ker(A^*$) for any square matrix $A$. This is not true in an infinite dimensional Hilbert space, in general. For example, if $A$ is the right shift operator on $\ell^2$, then $\dim\ker(A)=0\neq 1=\dim\ker(A^*)$.\\

\item Every matrix has a finite number of eigenvalues while an operator may have infinitely (even uncountably) many eigenvalues. For example, every $\lambda$ in the open unit disk of the complex plane is an eigenvalue of the left shift operator \cite[Example 2.3.2]{MUR}. On the other hand, the right shift operator has no eigenvalues. \\

\item Unlike the finite dimensional case in which the trace of each matrix is a complex number, the trace of an arbitrary operator $A\in  \mathbb{B}(\ell_2)$ defined by ${\rm tr}(A)=\sum_{j=1}^\infty\langle Ae_j,e_j\rangle$ may be infinite (or even non-existing). For example, for the diagonal operator \[A(x_1, x_2, x_3, \ldots)=(x_1, \frac{1}{2}x_2, \frac{1}{3}x_3, \ldots)\] on $\ell_2$, we have ${\rm tr}(A)=\sum_{j=1}^\infty \frac{1}{j}=\infty$. By the way, Grothiendieck \cite{GRO} has an example of an operator on a Banach Space where the trace is not the sum of the eigenvalues.\\

\item The spectrum of a matrix $A$ is contained in its numerical range, and the latter set is closed. Generally, neither statement is true for operators. For example, if $A$ is the diagonal operator ${\rm diag}(1, 1/2, 1/3, \ldots)$ (see the item above), then $\sigma(A)=\{1/n: n\in\mathbb{N}\}\cup\{0\} \nsubseteq (0,1]=W(A)$ and $W(A)$ is not a closed subset of the complex plane; cf. \cite[Problem 212]{HAL}. However, $\sigma(A)$ is a subset of the closure of $W(A)$ for every operator $A$.\\

\item Two operators $T$ and $S$ are similar if $T=W^{-1}SW$ for some invertible operator $W$. They are asymptotically similar if there exist sequences $(W_n)$ and $(V_n)$ of invertible operators such that $S=\lim_nW_n^{-1}TW_n$ and $T=\lim_nV_n^{-1}SV_n$. In the finite dimensional case, these two notions coincide but that is not the case in the infinite dimensional realm; cf. \cite[Theorem 2.1]{HER}. \\

\item It is known that the numerical range of any operator $A$ satisfying $A^n = I$ cannot be a disk in the finite dimensional setting; cf. \cite{MAR}. However, the authors of \cite{SPI1} construct an operator $A$ acting on an infinite dimensional Hilbert space such that $T^3 = I$ and $W(A)$ is an open disk centered at the origin.\\

\item If an invertible matrix $A$ is such that $||A^k\|,\,\, k=\pm1, \pm2, \ldots$ is constant, then $A$ is unitary. This is not so in the infinite-dimensional Hilbert spaces. Indeed, it is shown in \cite{GWW} that for each $\varepsilon>0$, there exists a nonunitary invertible operator $A$ on $\ell_2(\mathbb{Z})\oplus\ell_2(\mathbb{Z})$ such that $\|A^k\|= 1+\varepsilon$ for all $k\geq 1$.\\

\item The determinant of a matrix is equal to the product of its eigenvalues counted with their multiplicities. Evidently, this definition does not carry over to `all' operators acting on infinite dimensional Hilbert spaces.\\		
An extension of the notion of determinant is the Fredholm determinant, which is defined for operators of the form $I+A$, as an extension of $\det(I+A)=\exp({\rm tr}(\log(I+A)))$, where $A$ is a trace class operator, that is, an operator on a Hilbert space $\mathscr{H}$ such that $\sum_{e\in \mathcal{E}}\langle |A|e,e\rangle<\infty$, where $\mathcal{E}$ is an arbitrary orthonormal basis and $|A|$ stands for the positive square root of $A^*A$. Indeed, for operators in $I+{\rm trace~class}$ the determinant is the product of eigenvalues (this is usually stated in terms of the trace being the sum of eigenvalues for trace class operators and called Lidskii's theorem); see \cite{SIM}.

\item It is easily observed from \[\|Ax\|\leq \|A\|\,\|x\| \quad {\rm and} \quad \|Ax\|=\|\sum_{j=1}^n\langle x,e_j\rangle Ae_j\|\leq n\|x\|\max_{1\leq j\leq n}\|Ae_j\|\]
that a sequence $\{A_n\}$ converges to $A$ in the norm topology if and only if $\{A_nx\}$ converges to $Ax$ for all $x\in \mathbb{C}^n$. In infinite dimensional Hilbert spaces, the pointwise convergence does not imply the norm convergence, in general. For example, let $A_n\in  \mathbb{B}(\ell_2)$ be defined by the infinite diagonal matrix ${\rm diag}(1,1,\ldots, 1, 0,0,\ldots)$, whose first $n$ diagonal entries are equal to $1$, and all other entries are $0$. Then clearly $A_nx\to Ix$, for all $x\in \ell_2$, but the sequence $(A_n)$ is not a Cauchy sequence in  $\mathbb{B}(\ell_2)$ since $\|A_n-A_m\|=1,\,\, n\neq m$, and so cannot be convergent in the norm topology.\\

\item Every linear mapping  $A : \mathbb{C}^n \to \mathbb{C}^n$ is automatically continuous while a linear mapping on an infinite dimensional inner product space may be discontinuous (unbounded). Suppose that $\mathscr{K}$ is the dense subspace of $\ell_2$ consisting of all sequences $(x_n)$ with $x_n = 0$ for sufficiently large $n$. Let $A:\mathscr{K}\to \mathscr{K}$ denote the linear mapping $(x_n) \mapsto (nx_n)$. Then $A$ is unbounded since if $(e_n)$ is the orthonormal basis for $\ell_2$, then $\| e_n \| = 1$ and $\| Ae_n\| = n$ for all $n$.\\

In this example, $A$ is defined on a dense subset of $\ell_2$ but not on the whole space. Discontinuous linear operators defined on the whole space also exist and can be constructed with the use of Hamel bases. For example, following
Halmos \cite{HAL}: Extend the standard orthonormal basis $(e_n)$ of $\ell_2$ to a Hamel (linear algebra) basis $\beta$ for $\ell_2$. Choose $f \in \beta$ different from all $e_n$s, and define the linear operator $A: \ell_2\to \ell_2$ by $$A(g) = \left \{ \begin{array}{ll} 1& g=f\\ 0& g \in \beta\setminus\{f\} \end{array}\right .$$
Then $A(e_n) = 0$ and $A$ is unbounded (otherwise, $1=A(f)=\displaystyle{\sum_{n=1}^\infty}\langle f,e_n\rangle Ae_n = 0$).\\

\item Given an operator $A$, the unique operator $A^\dagger$ (if exists) satisfying (i) $AA^\dagger A=A$, (ii) $A^\dagger AA^\dagger=A^\dagger$, (iii) $A^\dagger A$ is self-adjoint, and (iv) $AA^\dagger$ is self-adjoint, is called the Moore--Penrose inverse of $A$. Every matrix has the Moore--Penrose inverse. However, there are operators having no Moore-Penrose inverses (precisely, those operators with non-closed ranges; see \cite{MSM1}). For example, the range of the operator $A$ on $\ell_2$ defined by $A(x_1, x_2, x_3, \ldots)=(x_1, \frac{1}{2}x_2, \frac{1}{3}x_3, \ldots)$ contains all finitely nonzero sequences, and so is dense in $\ell_2$. Since this range does not contain the sequence $(1/n)$, it is not closed. This $A$ has no Moore-Penrose inverse.\\

\item  It is known that every normal matrix can be written of the form $A=UDU^*$ with the unitary matrix $U$ and diagonal matrix $D={\rm diag}(\lambda_1, \ldots, \lambda_n)$, where $\lambda_1, \ldots, \lambda_n$ are the eigenvalues of $A$; see \cite[Theorem 9.1]{ZHA}. Such a result does not hold in the infinite dimensional case. In other words, there exist normal operators $A$ on an infinite dimensional Hilbert space $\mathscr{H}$ for which there are no orthonormal bases of $\mathscr{H}$ consisting of the eigenvectors of $A$. As an extreme manifestation of this phenomenon, the bilateral shift operator $A(f_n)=f_{n+1}\,\,(n=0, \pm 1, \pm 2, \ldots)$ on $\ell_2(\mathbb{Z})$ is normal but has no eigenvalues \cite[p. 56]{MUR}.\\

\item An operator $A$ is called hypercyclic if there exists a vector $x_0\in\mathscr{H}$ such that the set $\{A^nx_0: n=0, 1, 2, \ldots\}$ is dense in $\mathscr{H}$. The space
$\mathscr{H}$ must be separable in order to support a hypercyclic operator.
By using the Jordan decomposition for a matrix, it is not difficult to prove that if $\mathscr{H}$ is finite dimensional, then it has no hypercyclic operator
(see \cite[pp.~54--55]{GroP}). The situation for Hilbert spaces of infinite dimension is quite different. For example, every scalar multiple $\alpha A\,\, (|\alpha|>1)$ of the left shift operator $A$ on $\ell_2$ is a hypercyclic operator; see \cite{ROL}. A close, weaker notion is the one of supercyclicity. An operator $A$ is called hypercyclic if there exists a vector $x_0\in\mathscr{H}$ such that the projective orbit  $\{\lambda A^nx_0: n=0, 1, 2, \ldots; \lambda \in \K\}$ ($\K = \R$ or $\C$) is dense in $\mathscr{H}$.
Since any hypercyclic operator is supercyclic, we have that every infinite dimensional separable Hilbert space supports supercyclic operators (in fact, $A$ itself is supercyclic on $\ell_2$). In the finite dimensional case, Herzog \cite{HERZ} proved in 1992 that, if $\K = \R$ ($\K = \C$, resp.), then $\mathscr{H}$ supports a supercyclic operator if and only if ${\rm dim} (\mathscr{H}) \in \{1,2\}$ (${\rm dim} (\mathscr{H}) = 1$, resp.). In fact, every rotation on $\R^2$ given by a matrix
$\begin{pmatrix}
  \cos (2\pi \alpha ) &  \sin (2\pi \alpha )  \\
   -\sin (2\pi \alpha )  &  \cos (2\pi \alpha )
\end{pmatrix}$
with $\al$ irrational is supercyclic. \\

\item Factorization of matrices and operators acting on Hilbert spaces is a lively area of research in matrix analysis and operator theory. Problems of factorization ask whether a given operator in $\mathbb{B}(\mathscr{H})$ can be factored into (real or complex) linear combination or product of finitely many operators in a class of operators and seek for the minimal number of factors in a factorization. Matrix versions of these problems have a long history and many of them have appropriate analogues (probably under some additional conditions) for operators acting on Hilbert spaces of arbitrary dimension. However, some of these problems having solutions for matrices cannot have any solution for operators acting on infinite dimensional Hilbert space. A nice survey on such problems is \cite{WU}. By the polar decomposition, every matrix $A=U|A|$ is the product of two normal matrices, say $U$ and $|A|$, whilst the right shift operator cannot be factored as the product of finitely many normal operators; cf. \cite[Problem 144(a)]{HAL}.\\

\item  Bart et al. \cite{BAR} showed that if $P_1, \ldots, P_k$ are idempotent matrices such that $P_1+\ldots+P_k=0$, then $P_j=0$ for all $j=1, \ldots, k$. The situation changes in the infinite dimensional setting. As shown in \cite{BAR}, for $k\geq 5$ there exist $k$ different nonzero projections $P_1, \ldots, P_k$ on $\mathscr{H}$ such that $P_1+\ldots+P_k=0$. By the way, the number $5$ is sharp in the sense that there is no nontrivial zero sum of four idempotents.\\

\item For a long time, there has been considerable interest in the famous invariant subspace problem. This problem asks whether every operator $\mathfrak{T}$ on a Banach space $\mathscr{X}$ has a nontrivial (neither $\{0\}$ nor $\mathscr{X}$) invariant closed subspace. By an invariant subspace, we mean a subspace $\mathscr{M}$ such that $\mathfrak{T}\mathscr{M}\subseteq \mathscr{M}$. Enflo \cite{ENF} in 1975 proved that there exists a separable Banach space $\mathscr{X}$ and a continuous linear operator on $\mathscr{X}$ with dense range having no nontrivial closed invariant subspace. If $\mathscr{H}$ is a nonseparable Hilbert space, $x_0\neq 0$, and $A\in \mathbb{B}(\mathscr{H})$, then the closed linear span $\{A^nx_0: n=0, 1, 2, \ldots\}$ is a nontrivial invariant subspace  for $A$. By the spectral theorem, all normal operators on an infinite dimensional Hilbert space admit nontrivial invariant subspaces. The problem, in its generality, remains still open for (separable) Hilbert spaces. However, if $A\in \mathbb{M}_n\,\, n\geq 2$ is a matrix and $\lambda$ is an eigenvalue of $A$, then its eigenspace $\{x\in\mathbb{C}^n: Ax=\lambda x\}$ is a nontrivial invariant subspace for $A$.
\end{itemize}
\smallskip
\textbf{Acknowledgement.} The authors would like to, sincerely, thank Professor Ilya M. Spitkovsky for his valuable comments improving this note.


\begin{thebibliography}{99}

\bibitem{BAR} H. Bart, T. Ehrhardt, and B. Silbermann, \textit{Zero sums of idempotents in Banach algebras}, Integral Equations Operator Theory \textbf{19} (1994), no. 2, 125--134.

\bibitem{BHA} R. Bhatia, \textit{Positive definite matrices}, Princeton Series in Applied Mathematics. Princeton University Press, Princeton, NJ, 2007.

\bibitem{ENF} P. Enflo, \textit{On the invariant subspace problem for Banach spaces}, Acta Math. \textbf{158} (1987), no. 3-4, 213--313.

\bibitem{FRR} C. K. Fong, H. Radjavi, and P. Rosenthal, \textit{Norms for matrices and operators}, J. Operator Theory \textbf{18} (1987), no. 1, 99--113.

\bibitem{GWW} H.-L. Gau, K.-Z. Wang, and P. Y. Wu, \textit{Constant norms and numerical radii of matrix powers}, Oper. Matrices \textbf{13} (2019), no. 4, 1035--1054.

\bibitem{GO} B. R. Gelbaum and J. M. H. Olmsted, \textit{Theorems and counterexamples in mathematics}, Problem Books in Mathematics. Springer-Verlag, New York, 1990.

\bibitem{GroP} K. G. Grosse-Erdmann and A. Peris, \textit{Linear Chaos}, Springer, London, 2011.

\bibitem{GRO} A. Grothendieck, \textit{Produits tensoriels topologiques et espaces nucl\'{e}aires (French)}, Mem. Amer. Math. Soc. 16 (1955), Chapter 1: 196 pp.; Chapter 2: 140 pp.

\bibitem{SPI1} T. R. Harris, M. Mazzella, L. J. Patton, D. Renfrew, and I. M. Spitkovsky, \textit{Numerical ranges of cube roots of the identity}, Linear Algebra Appl. \textbf{435} (2011), no. 11, 2639--2657.

\bibitem{HAL} P. R. Halmos, \textit{A Hilbert space problem book}, Second edition. Graduate Texts in Mathematics, 19. Encyclopedia of Mathematics and its Applications, 17. Springer-Verlag, New York-Berlin, 1982.

\bibitem{HER} D. A. Herrero, \textit{Approximation of Hilbert space operators}, Vol. 1. Second edition. Pitman Research Notes in Mathematics Series, 224. Longman Scientific \& Technical, Harlow; copublished in the United States with John Wiley \& Sons, Inc., New York, 1989.

\bibitem{HERZ} G. Herzog, \textit{On linear operators having supercyclic vectors}, Studia Math. \textbf{103} (1992), no. 3, 295--298.

\bibitem{MAR} J. Maroulas, P. J.  Psarrakos, and M. J. Tsatsomeros, \textit{Perron-Frobenius type results on the numerical range}, Linear Algebra Appl. \textbf{348} (2002), 49--62.

\bibitem{MOS2} M. S. Moslehian, \textit{Trick with $2\times2$ matrices over $C^*$-algebras}, Austral. Math. Soc. Gaz. \textbf{30} (2003), no. 3, 150--157.

\bibitem{MSM1} M. S. Moslehian, K. Sharifi, M. Forough, and M. Chakoshi, \textit{Moore-Penrose inverses of Gram operators on Hilbert $C^*$-modules}, Studia Math. \textbf{210} (2012), no. 2, 189--196.

\bibitem{MUR} G. J. Murphy, \textit{$C^*$-algebras and operator theory}, Academic Press, INC, 1990.

\bibitem{ROL} S. Rolewicz, \textit{On orbits of elements}, Studia Math. \textbf{32}  (1969), 17--22.

\bibitem{SIM} B. Simon, \textit{Operator theory. A Comprehensive Course in Analysis}, Part 4. American Mathematical Society, Providence, RI, 2015.

\bibitem{WU} P. Y. Wu, \textit{The operator factorization problems}, Linear Algebra Appl. \textbf{117} (1989), 35--63.

\bibitem{ZHA} F. Zhang, \textit{Matrix theory. Basic results and techniques}, Second edition. Universitext. Springer, New York, 2011.

\end{thebibliography}
\end{document}